\documentclass{article}

\usepackage{arxiv}

\usepackage[utf8]{inputenc} 
\usepackage[T1]{fontenc}    
\usepackage{hyperref}       
\usepackage{url}            
\usepackage{booktabs}       
\usepackage{amsfonts}       
\usepackage{nicefrac}       
\usepackage{microtype}      
\usepackage{lipsum}		
\usepackage{graphicx}
\usepackage[square,sort,comma,numbers]{natbib}
\usepackage{doi}

\usepackage{amssymb}
\usepackage{listings}
\usepackage{amsmath}
\usepackage{paralist}
\usepackage[utf8]{inputenc}
\usepackage[english]{babel}

\usepackage{subfig}
\usepackage{epsfig} 
\usepackage{epstopdf}

\newcommand{\norm}[1]{\left\lVert#1\right\rVert}

\title{Smart Gradient - An Adaptive Technique for Improving Gradient Estimation}

\author{ \href{https://orcid.org/0000-0003-1587-3288}{\includegraphics[scale=0.06]{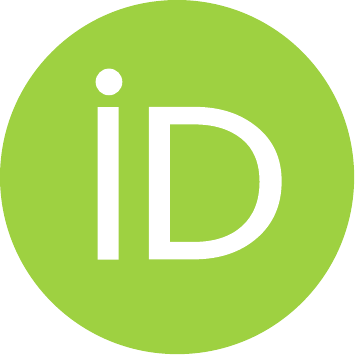}\hspace{1mm}Esmail H.~Abdul Fattah}\thanks{Use footnote for providing further
		information about author (webpage, alternative
		address)---\emph{not} for acknowledging funding agencies.} \\
	Department of Statistics\\
	King Abdullah University of Science and Technology\\
	Thuwal, 23955, Makkah\\
	\texttt{esmail.abdulfattah@kaust.edu.sa} \\
	\And
	\href{https://orcid.org/0000-0002-4334-2057}Janet ~Van Niekerk \\
	Department of Statistics\\
	King Abdullah University of Science and Technology\\
	Thuwal, 23955, Makkah \\
	\texttt{janet.vanniekerk@kaust.edu.sa} \\
	\And
	\href{https://orcid.org/0000-0002-0222-1881}H{\aa}vard ~ Rue \\
	Department of Statistics\\
	King Abdullah University of Science and Technology\\
	Thuwal, 23955, Makkah\\
	\texttt{haavard.rue@kaust.edu.sa} \\

}

\renewcommand{\shorttitle}{Smart Gradient}

\begin{document}
\maketitle

\begin{abstract}

Computing the gradient of a function provides fundamental information about its behavior. This information is essential for several applications and algorithms across various fields. One common application that require gradients are optimization techniques such as stochastic gradient descent, Newton's method and trust region methods. However, these methods usually requires a numerical computation of the gradient at every iteration of the method which is prone to numerical errors. We propose a simple limited-memory technique for improving the accuracy of a numerically computed gradient in this gradient-based optimization framework by exploiting (1) a coordinate transformation of the gradient and (2) the history of previously taken descent directions. The method is verified empirically by extensive experimentation on both test functions and on real data applications. The proposed method is implemented in the \textbf{\texttt{R}} package \textbf{\texttt{smartGrad}} and in C\texttt{++}.

\end{abstract}


\section{Introduction}\label{Introduction_section}

Gradients are one of the oldest constructs of modern mathematics and often form the basis for many optimization problems. Moreover, gradients are used in various function expansions like the Taylor series expansion, in tangent plane construction and in various real-life engineering challenges such as building ramps, roads and buildings, amongst others. The gradient of a mathematical function $f$, at the input $\pmb{x}$, is denoted by $\nabla f(\pmb{x})$. It can be mathematically interpreted as a rate of disposition based on the disposition of $\pmb{x}$ or graphically as the slope of the tangent plane at $\pmb{x}$ \cite{Nocedal1999NumericalO}. Often though, the analytical form of $\nabla f(\pmb{x})$ is unknown or computationally intensive to evaluate and hence the popularity of numerical gradients methods. Newton's and quasi-Newton methods \cite{Fletcher1988PracticalMO} are well-known frameworks in optimization that use (numerical) gradients to get the descent directions.

Consider the problem of minimizing a twice differentiable continuous function $f$,
where $f(\pmb x): \mathbb{R}^n \xrightarrow{} \mathbb{R}$, is convex, i.e. satisfies,
$$f((1 - \lambda) \pmb x + \lambda \pmb y) \leq (1 - \lambda) f(\pmb x) + \lambda f(\pmb y) $$
for $\pmb x$, $\pmb y$ $\in$ $\mathbb{R}^n$ and $\lambda \in [0,1]$. Assume this (unconstrained) optimization problem is solvable, where its optimum $\pmb x^*$ exists and is unique. In this scenario, a necessary and a sufficient condition for $\pmb x^*$ to be optimal is
$$\nabla f(\pmb x^{*}) = \pmb 0.$$ 
\noindent Here, $\nabla f(\pmb x)$ is composed of the function's partial derivatives, that describe the rates of change in multiple dimensions. More generally we consider the directional derivative of a function. Given a continuous function $g$ with its first order partial derivatives, the directional derivative \citep{Thomas2005ThomasCE} of $g$ at $\pmb x$ $\in$ $\mathbb{R}^n$ in the direction of unit vector $\pmb u$ is
$$D_{\pmb u}g(\pmb x) = \frac{\partial g}{\partial x_1} u_1 + \frac{\partial g}{\partial x_2} u_2 + \ldots + \frac{\partial g}{\partial x_n} u_n.$$

The gradient vector can be reformulated as a combination of directional derivatives based on the canonical basis $\{\pmb e_1, \ldots, \pmb e_n \}$, where every element of $\pmb e_i$ is zero except for the $i^{\text{th}}$ position being 1, as in formula \eqref{gradform1}, that will follow shortly.

Common methods for solving optimizations problems rely on iterative techniques which generate iterations that converge to the optimum $\pmb x^*$. The iterations are constructed from the general iterative update equation,
$$\pmb x^{(k+1)} = \pmb x^{(k)} + \alpha_k \pmb d^{(k)},$$ with step size $\alpha_k$ in the direction of the vector $\pmb d^{(k)}$. The choice of the directions $\{\pmb d_k\}_{k=1}^{n}$ depends on the utilized method \cite{Nocedal2020NumericalO2}:
\begin{enumerate}
    \item Gradient Descent: $\pmb d^{(k)} = -\nabla f(\pmb x^{(k)})$.
    \item Newton's Method: $\pmb d^{(k)} = - \nabla^2 f( \pmb x^{(k)})^{-1} \nabla f(\pmb x^{(k)})$.
    \item Quasi Newton's Method: $\pmb d^{(k)} = - \pmb B_k \nabla f(\pmb x^{(k)})$ where $\pmb B_k$ is an estimate of the Hessian at $\pmb x^{(k)}$, $\nabla^2 f(\pmb x^{(k)})$.
\end{enumerate}
The choice of the step size $\{\alpha^{(k)}\}$ is obtained using either exact or inexact line search \cite{Nocedal1999NumericalO}.

We summarize and relate four popular gradient computational frameworks next \cite{depner2016hydrodynamics}.

\begin{enumerate}
    \item \textbf{Exact gradient with the canonical basis:} Based on directional derivatives, the gradient can be computed using the canonical basis $\{\pmb e_1, \ldots, \pmb e_n \}$,
    \begin{equation}
 	    \nabla f(\pmb x) = \Big(D_{\pmb e_1}f(\pmb x), D_{\pmb e_2}f(\pmb x) , \ldots, D_{\pmb e_n}f(\pmb x)\Big).
 	\label{gradform1}
    \end{equation}
    
    \item \textbf{Exact gradient with a non-canonical basis:}
    Here, the gradient is computed based on the directional derivatives using a set of directions obtained from $\{\pmb v_{1}, \ldots, \pmb v_{n}\}$. For a non-singular matrix $\pmb G = [\pmb v_{1}|\pmb v_{2}|\ldots|\pmb v_{n}]$, then 
    \begin{equation}
        \nabla_{\pmb v} f(\pmb x) = {\pmb G}^{-T} \nabla h(\pmb \varphi)\Big|_{\pmb \varphi = 0} \text{ where } h(\pmb \varphi) = f(\pmb  x + \pmb G \pmb \varphi).
        \label{gradform2}
    \end{equation}

    \item \textbf{Inexact gradient in canonical basis: Vanilla Gradient (VG)}. The gradient of the objective function is computed numerically as an estimate to the exact gradient with the canonical basis, such as using finite-difference methods,
    \begin{equation}
        \Tilde{\nabla} f(\pmb x) \approx {\nabla} f(\pmb x).
        \label{gradform3}
    \end{equation}

    \item \textbf{Inexact gradient in non-canonical basis:} The gradient is computed using a general basis instead of the canonical basis, 
        \begin{equation}
        \Tilde{\nabla}_{\pmb v} f(\pmb x) = {\pmb G}^{-T} \Tilde{\nabla} h(\pmb \varphi)\Big|_{\pmb \varphi = 0} \approx {\nabla} f(\pmb x) \text{ where } h(\pmb \varphi) = f(\pmb  x + \pmb G \pmb \varphi).
        \label{gradform4}
    \end{equation}
\end{enumerate}

\noindent For $\pmb G^T = \pmb I_n$, \eqref{gradform2} reduces to \eqref{gradform1}, and \eqref{gradform4} reduces to the Vanilla Gradient in \eqref{gradform3} where the partial derivatives are directly computed along the vectors of the canonical basis.

\begin{figure}[!htb]
    \centering
    \includegraphics[scale=0.2]{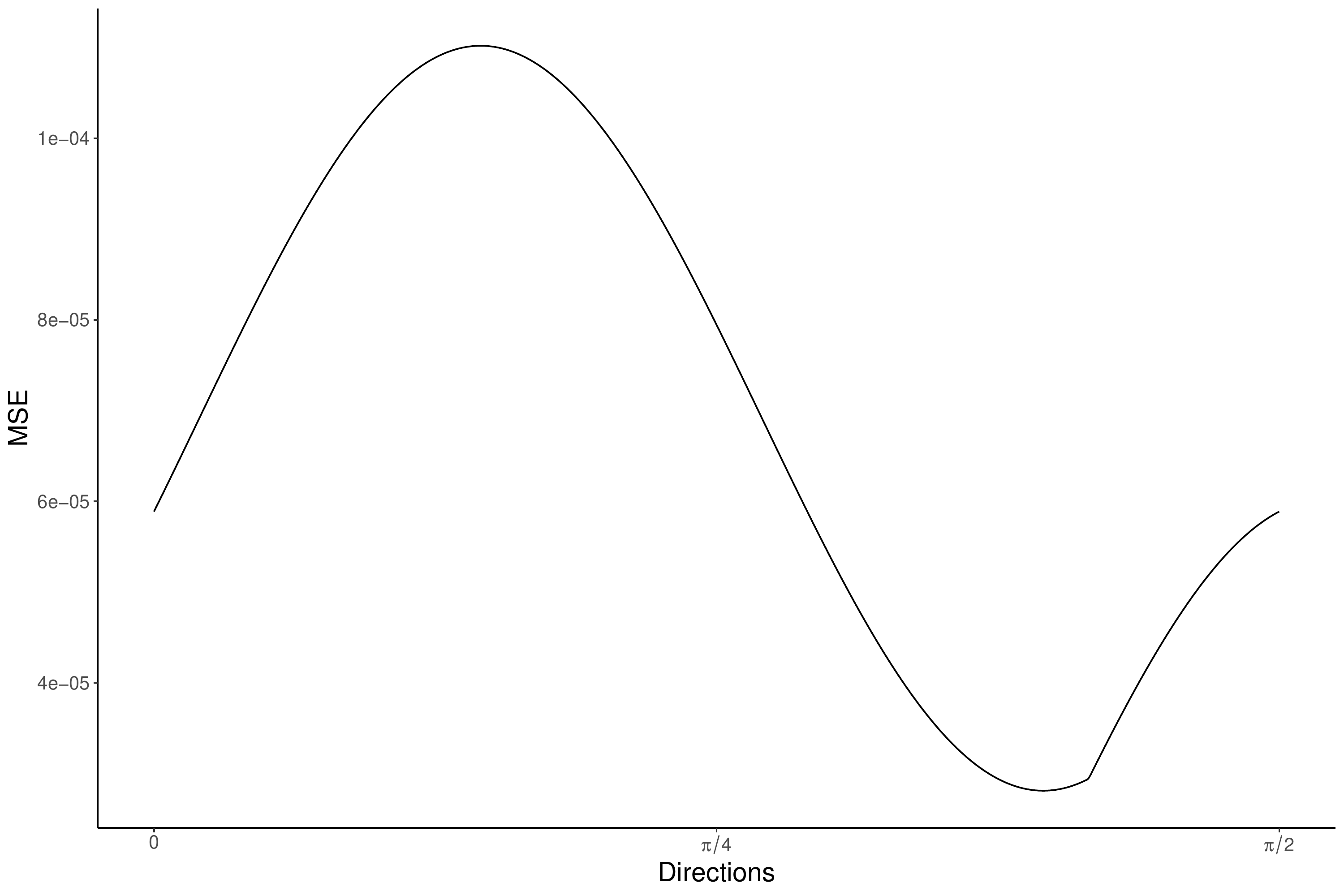}
    \includegraphics[scale=0.2]{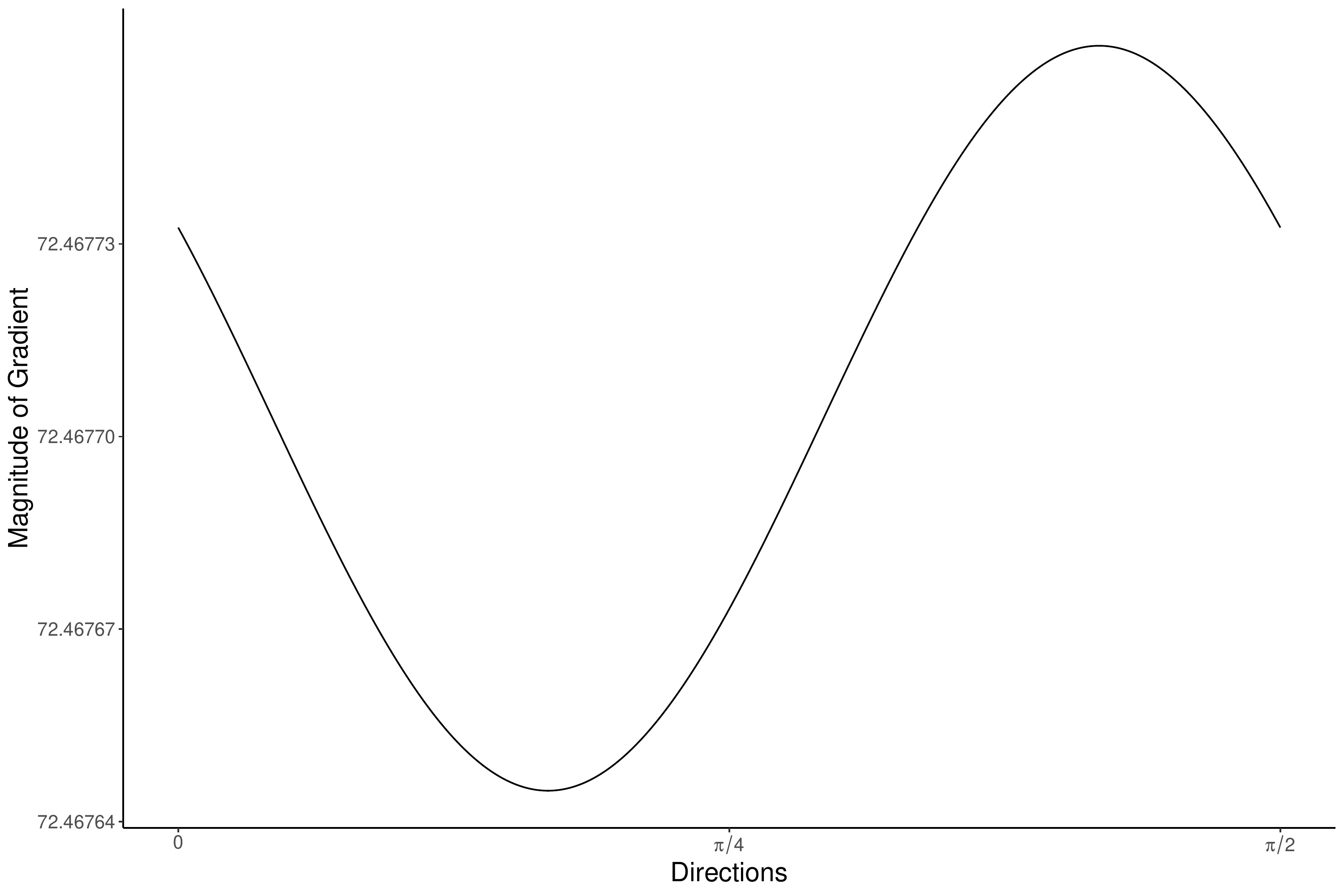}
    \caption{Estimating Gradient using different Directions}
    \label{fig:basis_rotation}    
\end{figure}

A natural question is if we can construct $\pmb G$ in \eqref{gradform4} that results in a more accurate Vanilla Gradient. As an illustration, we estimate the gradient of the two-dimensional Rosenbrock function \cite{besag1975statistical}
$f(\pmb x) = (1 - x_1)^2 + 100 \times (x_2 - x_1^2)^2$, see Figure \ref{fig:RosenFun}, using various bases. We use the central difference method of first order of step length $10^{-3}$ to get the gradient estimate at $\pmb x^{(0)} = (-0.29, 0.40)^T $. As for the directions, we start with the unit vectors $\pmb e_1 = (1,0)^T$ and $\pmb e_2 = (0,1)^T$ then we keep rotating with an angle of $\pi/10^3$ until most of the directions are explored. The mean square error for the estimated gradient is calculated for each set of directions, and at each direction the magnitude of the gradient is calculated. The results are presented in Figure \ref{fig:basis_rotation}.

From Figure 1 we observe opposite trends between the error in the
gradient estimates and the magnitude of the gradient itself: the MSE
error in the gradient estimate is low when the magnitude of the gradient
is high. Although such a claim does not hold in general, we found
empirically it holds ``almost everywhere", and this is sufficient
for our proposed method. This observation, suggests the obvious
approach, here explained in dimension 2 for simplicity, to improve the
gradient estimates at $\pmb x$,
\begin{enumerate} 
    \item Find the direction $\pmb d$ that maximize the change in the function itself. 
    \item Estimate the gradient in direction $\pmb d$ and the direction orthogonal to $\pmb d$. 
    \item Do a linear transformation for the estimates we get in 2 to obtain the estimates of the Vanilla Gradient as in \eqref{gradform4}. \end{enumerate}

We are moving around the same point: to get the descent direction we estimate the gradient and to estimate the gradient we use direction where there is a high change in the objective value $f$, and this reduction is exactly what we want for the descent direction. This means that the  directions  we  should  use  for  estimating  the  gradient  are  exactly  what  we want to obtain by computing gradients. Within an iterative optimization
algorithm, we have access to $\pmb x$ at most recent iterations,
$\pmb x^{(k)},\pmb  x^{(k-1)}, \ldots$, then the most recent change in $\pmb x$ can
be used as the (surrugate) direction $\pmb d^{(k)}= \pmb x^{(k)}- \pmb x^{(k-1)}$. In
dimension $n$, we would need to use the $n$ most recent differences as
surrugate directions. These directions should still be relevant, and hence we assume a low or moderate dimension of $\pmb x$.  \\

In this paper, we propose in Section \ref{Methodology_section} a method for constructing the matrix $\pmb G$ adaptively . In Section \ref{Numerical_section}, we demonstrate how this construction results in improving the accuracy of a numerical gradient using test functions and some applications within the R-INLA package. The paper is concluded by a discussion in Section \ref{Discussion_section}.

\section{Methodology}\label{Methodology_section}
The proposed method is based on using the prior information we have about the descent directions in previous iterations, $\{\pmb x^{(i)}\}_{i=k-n}^{k}$ for improving the estimated numerical gradient $\Tilde{\nabla} f(\pmb x)$. The difference between $\pmb x^{(i)}$ and $\pmb x^{(i-1)}$, $i = k, \ldots k-n$, indicates the direction where the objective function is highly reduced at that position. We then use these $n$ differences as (surrogate) directions to estimate the gradient.

At iteration $k$, the $n$ most recent differences can be used and the $n$ directions have the following form,
$$\pmb d^{(k)} = \frac{\Delta \pmb x^{(k)}}{\left \|  \Delta \pmb x^{(k)} \right \|}$$
where $\Delta \pmb x^{(k)} = \pmb x^{(k)} - \pmb x^{(k-1)}$. To get a unique contribution from each direction $\pmb d^{(k)}$, we subtract all the components of $\pmb d^{(k)}$ that are in the same direction of $\{\pmb d^{(j)}\}_{j=k-n+1}^{k-1}$, then normalize it. We can do this by expressing these directions in terms of projectors, 
$$\pmb {\Tilde{d}}^{(k)} = \pmb P_{\perp \pmb {\Tilde{d}}^{(k-1)}} \pmb P_{\perp \pmb {\Tilde{d}}^{(k-2)}} \dots \pmb P_{\perp \pmb {\Tilde{d}}^{(k-n+1)}} \pmb d^{(k)}$$
where $\pmb P_{\perp \pmb {\Tilde{d}}^{(k)}} = I - \pmb {\Tilde{d}}^{(k)}  (\pmb {\Tilde{d}}^{(k)})^T$. We apply this orthogonalization to the $n$ given directions using the Modified Gram-Schmidt (MGS) orthogonalization \cite{Picheny2013ABO}, which is an algorithm used to compute an orthonormal basis $\{\pmb {\Tilde{d}}_1^{(k)},\pmb {\Tilde{d}}_2^{(k)},\ldots,\pmb {\Tilde{d}}_n^{(k)}\}$ that spans the same subspace as the original vectors, i.e,
$$\texttt{Span}(\{\pmb d_1^{(k)},\pmb d_2^{(k)},\ldots,\pmb d_n^{(k)}\}) = \texttt{Span}(\{\pmb {\Tilde{d}}_1^{(k)},\pmb {\Tilde{d}}_2^{(k)},\ldots,\pmb {\Tilde{d}}_n^{(k)}\}).$$

\noindent The matrix $\pmb {\Tilde{G}}^{(k)}$ represents the orthogonal directions at iteration k,
$$\pmb {\Tilde{G}}^{(k)} = \Big[{ \Tilde{\pmb d}}_{1}^{(k)} | {\Tilde{\pmb d}}_{2}^{(k)} | \ldots | {\Tilde{\pmb d}}_{n}^{(k)} \Big],$$
and is the MGS of matrix $\pmb G^{(k)}$,
$$\pmb G^{(k)} = \Big[\Delta \pmb x^{(k)} | \pmb d_1^{(k-1)} | \pmb d_2^{(k-1)}  | \ldots | \pmb d_{n-1}^{(k-1)} \Big].$$
We use matrix $\pmb {\Tilde{G}}^{(k)}$ to estimate the gradient at $\pmb x^{(k)}$, using \eqref{gradform4},
 \begin{equation}
        \Tilde{\nabla}_{{\Tilde{\pmb d}}} f^{(k)}(\pmb x^{(k)}) = {{\pmb {\Tilde{G}}^{(k)}}}^{-T} \Tilde{\nabla}_{\pmb e} h^{(k)} (\pmb \varphi)\Big|_{\pmb \varphi = 0} \text{ where } h(\pmb \varphi) = f(\pmb  x^{(k)} + \pmb {\Tilde{G}}^{(k)} \pmb \varphi).
        \label{gradform5}
\end{equation}
We start with $\pmb {\Tilde{G}}^{(0)} = \pmb I_n$. This transformation \eqref{gradform5} is a rotation/reflection using a different basis, and the length of the gradient is invariant due to the property of orthogonal matrices. The simplicity of this technique is represented by its easy implementation. The gradient is just computed in a different coordinate system using a simple reparamaterization of function $f$ and a matrix-vector multiplication. We label this new approach by Smart Gradient. \\

\textbf{Simple example for illustration}
Assume we have two iterations of resulting in the set of $\pmb x$ positions of a two-dimensional function,
\begin{center}
$\pmb x^{(0)} = (1.78, 2.82)^T$, $\pmb x^{(1)} = (1.89, 4.62)^T$, $\pmb x^{(2)} = (11.54, 4.15)^T$, $\ldots$
\end{center}

At each position $k$, we form $\pmb G^{(k)}$. This matrix is initialized as the identity matrix, so at the first position $k = 0$ the first direction used is $\pmb d_1^{(0)} = (1,0)^T$ and the second $\pmb d_2^{(0)} = (0,1)^T$. When reaching the second position, we get a new direction $\Delta{\pmb x^{(1)}} = \pmb x^{(1)} - \pmb x^{(0)}$, so matrix $\pmb {G}^{(1)}$ is constructed as follow:
$$\pmb G^{(1)} = \Big[\pmb x^{(1)} - \pmb x^{(0)} |\pmb d_1^{(0)} \Big] = \begin{pmatrix} 0.11 & 1\\ 1.80 & 0 \end{pmatrix}$$

This matrix is orthonormalized using the MGS algorithm, such that the columns of $\pmb G^{(1)}$, $\pmb g_1^{(1)}$ and $\pmb g_2^{(1)}$, are updated to form $\pmb {\Tilde{G}}^{(1)}$ by columns ${\Tilde{\pmb d}}_{1}^{(1)}$ and ${\Tilde{\pmb d}}_{2}^{(1)}$,

\begin{center}
	${\Tilde{\pmb d}}_{1}^{(k)} = \dfrac{ \Delta{\pmb x^{(k)}} }{\norm{\Delta{\pmb x^{(k)}}}}$ and ${\Tilde{\pmb d}}_{2}^{(k)} = \dfrac{\Delta{\pmb x^{(k-1)}} - \Big( \pmb g_1^{(k)} \cdot \Delta{\pmb x^{k-1}}\Big) \dfrac{ \Delta{\pmb x^{(k-1)}} }{\norm{\Delta{\pmb x^{(k-1)}}}}  }{ \norm{\Delta{\pmb x^{(k-1)}} - \Big( \pmb g_1^{(k)} \cdot \Delta{\pmb x^{(k-1)}}\Big) \dfrac{ \Delta{\pmb x^{(k-1)}} }{\norm{\Delta{\pmb x^{(k-1)}}}}}},$
\end{center}
where $\pmb g_i^{(k)}$ is the ith column of $\pmb G^{(k)}$, and
\begin{center}
$
\pmb G^{(1)} = \begin{pmatrix}
0.0601 & 0.9981\\
0.9981 & -0.0610
\end{pmatrix}.
$
\end{center}
At the next position, $\pmb x_2$, the same procedure for the two updates are computed,
$$\pmb G^{(2)} = \Big[\pmb x^{(2)} - \pmb x^{(1)} |\pmb d_1^{(1)} \Big] = \begin{pmatrix} 9.65 & 0.0610\\ -0.47 & 0.9981 \end{pmatrix}$$
then after MGS, 
\begin{center}  
$\pmb{ \Tilde{G}}^{(2)} = 
\begin{pmatrix}
0.9989 & 0.0486\\
-0.0486 & 0.9989
\end{pmatrix}
$
\end{center}

The columns of $\pmb{ \Tilde{G}}^{(2)}$ are orthogonal unit vectors, and they are used as the new set of axes in an orthogonal coordinate system in this unitary transformation of the gradient. For each formed $\pmb{ \Tilde{G}}^{(k)}$, the gradient is estimated using \eqref{gradform5}.


\section{Numerical Experiments}\label{Numerical_section}

\subsection{Test Functions}
In this section, some experiments are carried out to compare the performance of using the proposed approach. We chose to compare our approach with another numerical gradient estimator based on the first order central differences.
\begin{figure}[hbt!]
    \centering
    \includegraphics[scale=0.16]{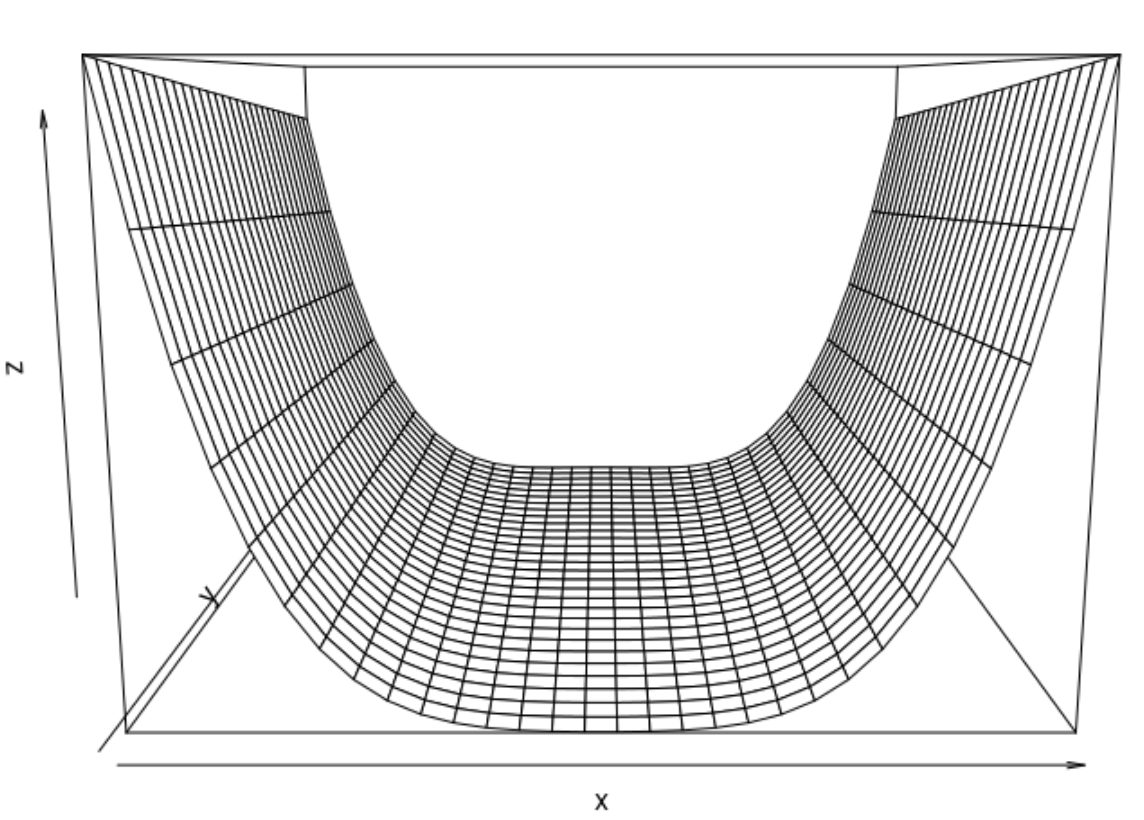}
    \caption{Rosenbrock Function}
    \label{fig:RosenFun}
\end{figure}

For each test function, we use \texttt{optim} function in \texttt{stats} package in R with the BFGS method \cite{Fletcher1988PracticalMO} to get the optimum, and the average squared difference (MSE) between the exact gradient (calculated analytically) and the estimated gradient (using Vanilla Gradient and the Smart Gradient) is calculated at each iteration. The experiment is repeated $100$ times with a random initial value.

One of the popular test problems for unconstrained optimization is the Extended Rosenbrock Function, which is unimodal, differentiable and has an optimum in a narrow parabolic valley, see Figure \ref{fig:RosenFun}. Another one is the Extended Freudenstein Roth Function, and the equations for both functions are here:
\begin{enumerate}
    \item Extended Rosenbrock Function:
    $$f(\pmb x) = \sum_{i=1}^{n/2} 100(x_{2i} - x_{2i-1}^2)^2 + (1-x_{2i-1})^2$$
    \item Extended Roth Freudenstein Function:
    \begin{equation}
        \begin{split}
            f(\pmb x) &= \sum_{i=1}^{n/2} \Big(-13 + x_{2i-1} + x_{2i}(x_{2i}(5 - x_{2i}) - 2)\Big)^2 
            \\&+ \Big(-29 + x_{2i-1} + x_{2i}(x_{2i}(x_{2i} +1 )  -14   ) \Big)^2
        \end{split}
    \end{equation}

\end{enumerate}

Comparison results are summarized in Table (\ref{Tab:RosenRoth}) showing an improvement (ratio of the errors of VG and SG)  of 2.5 for dimension 5 and at least 3.5 for the other two higher dimensions in Rosenbrock function. 

\begin{table}[]
\centering
\begin{tabular}{|c|c|c|c|}
\hline
\textbf{$\pmb x$ dimension} & \textbf{\begin{tabular}[c]{@{}c@{}}Average MSE \\ Vanilla Gradient\end{tabular}} & \textbf{\begin{tabular}[c]{@{}c@{}}Average MSE \\ Smart Gradient\end{tabular}} & \textbf{Improvement} \\ \hline
\multicolumn{4}{|c|}{\begin{tabular}[c]{@{}c@{}} \\ \textbf{Extended Rosenbrock Function} \end{tabular}}
\\ \hline
5             & 2.60e-04                                                                                 & 1.04e-04                                                                                & 2.5         \\ \hline
10            & 2.71e-04                                                                                 & 0.78e-04                                                                                & 3.47        \\ \hline
25            & 2.80e-04                                                                                 & 0.49e-04                                                                                & 5.71        \\ \hline
\multicolumn{4}{|c|}{\begin{tabular}[c]{@{}c@{}} \\ \textbf{Extended Freudenstein Roth Function} \end{tabular}}                                                                                                                                                        \\ \hline
5             & 2.26e-04                                                                                 & 1.39e-04                                                                                & 1.63        \\ \hline
10            & 2.78e-04                                                                                 & 1.42e-04                                                                                & 1.96        \\ \hline
25            & 2.84e-04                                                                                 & 1.25e-04                                                                                & 2.27        \\ \hline
\end{tabular}
\captionof{table}{The average MSE when using Vanilla and Smart Gradient approaches compared to the exact gradient for different dimensional Rosenbrock and Roth functions  \label{Tab:RosenRoth}}
\end{table}

\begin{figure}[hbt!]
    \centering
    \includegraphics[scale=0.23]{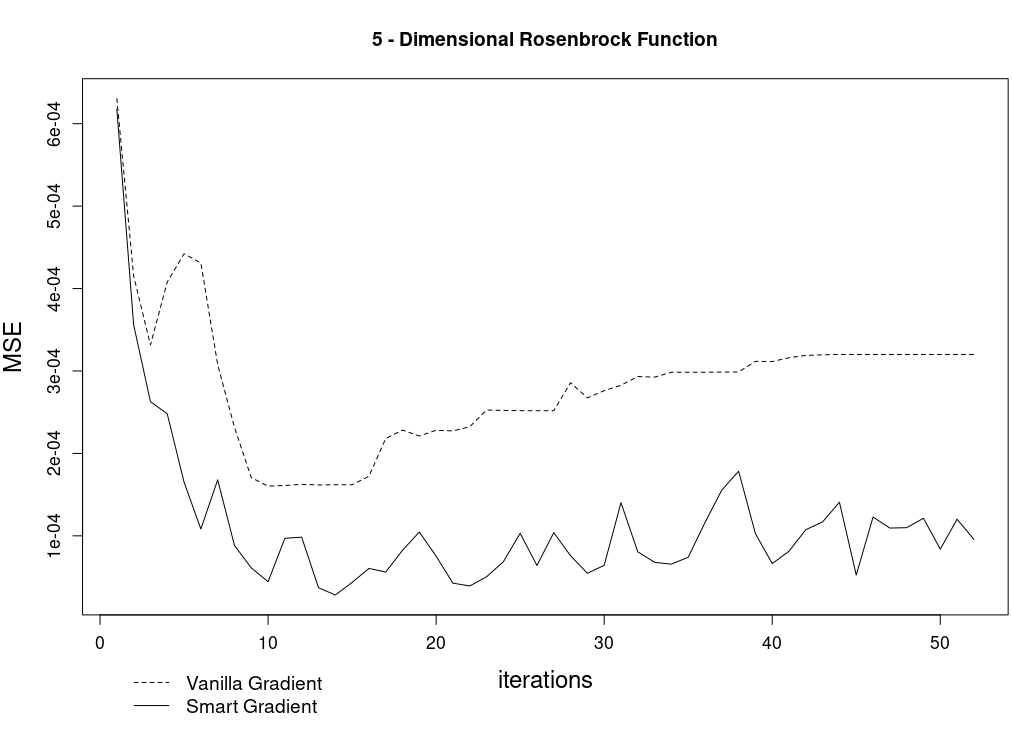}
    \includegraphics[scale=0.23]{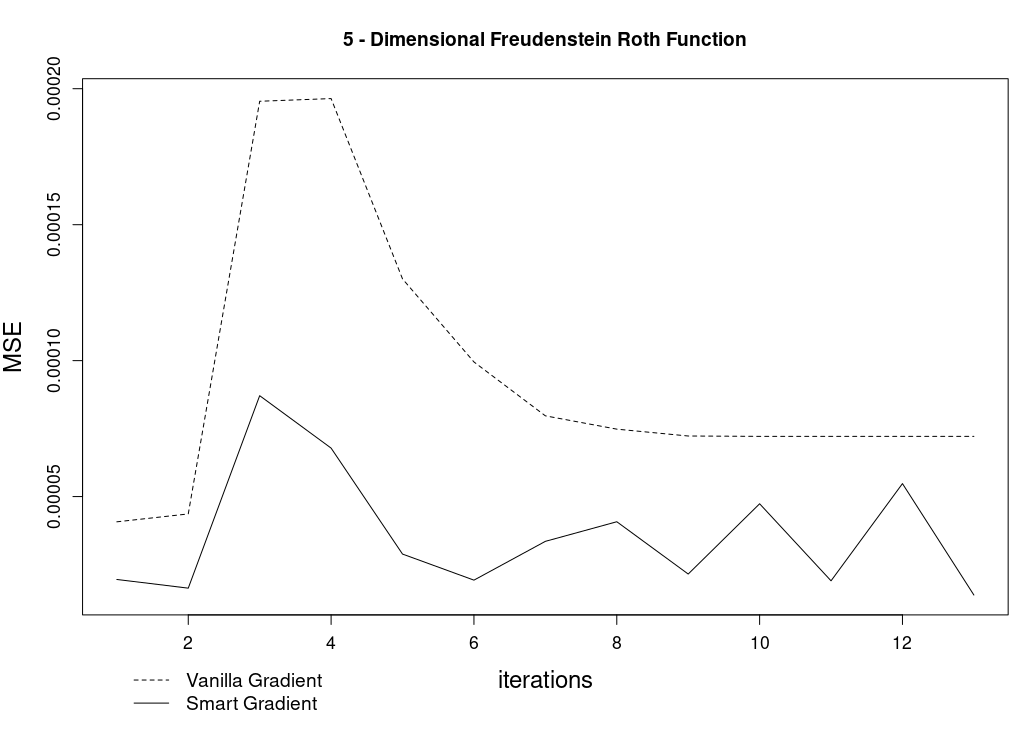}
    \includegraphics[scale=0.23]{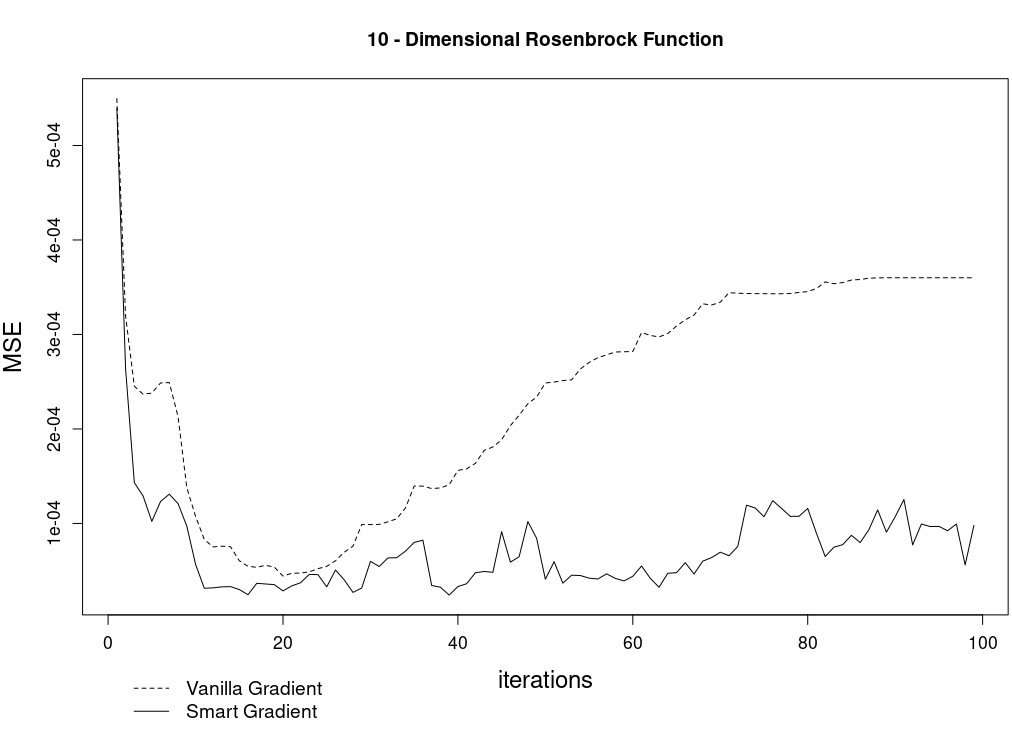}
      \includegraphics[scale=0.23]{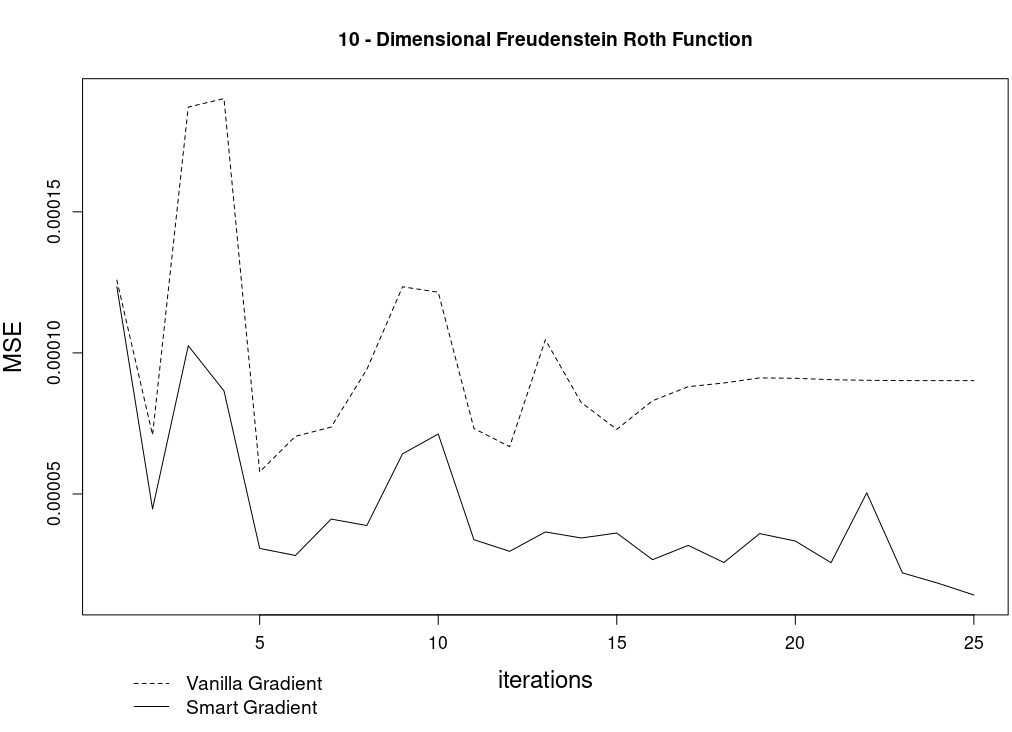}
    \includegraphics[scale=0.23]{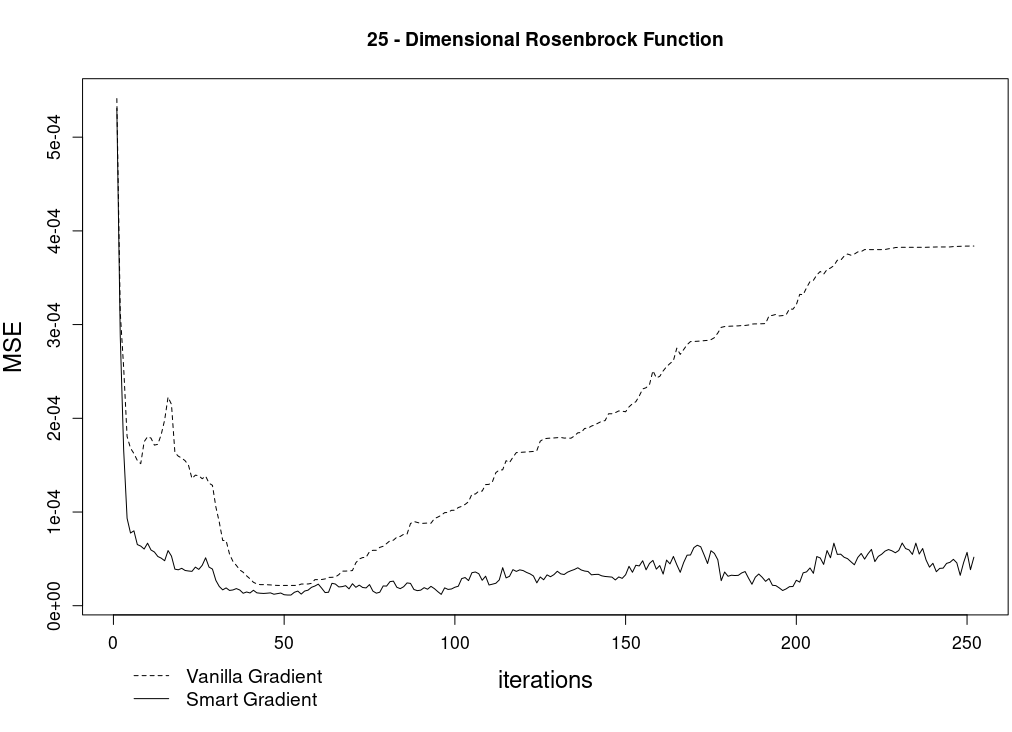}
    \includegraphics[scale=0.23]{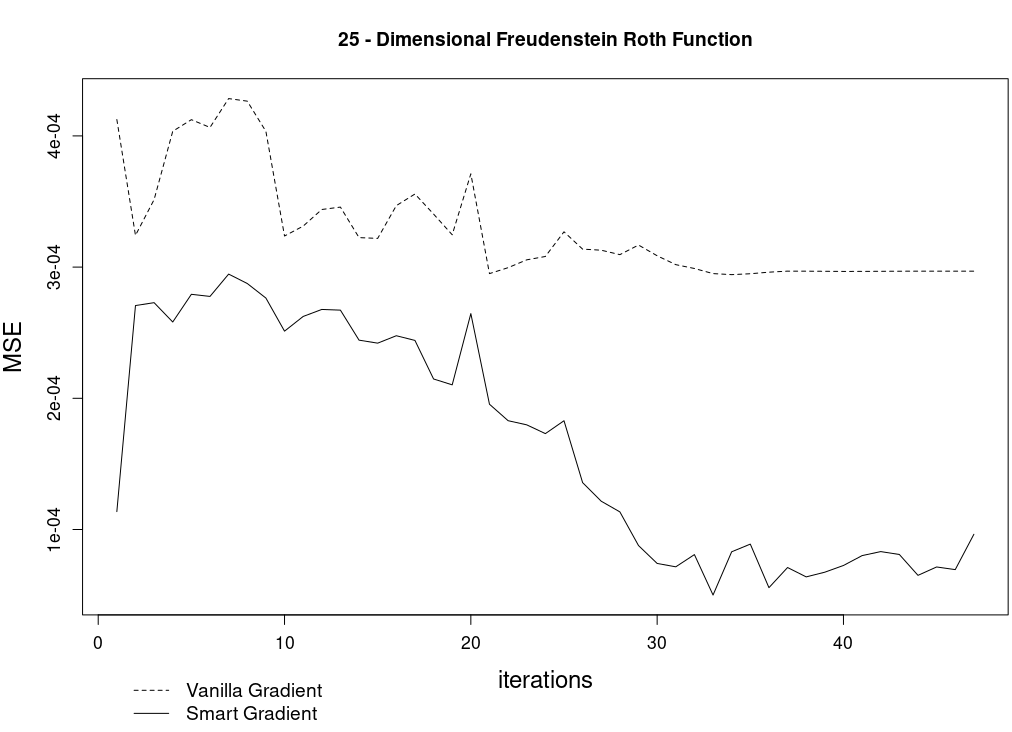}
    \caption{MSE for Smart and Vanilla Gradients at each iteration using Different Dimensional Rosenbrock and Roth Functions}
    \label{fig:RothandRothFun}
\end{figure}

In Figure \ref{fig:RothandRothFun}, the two estimated gradients started with almost the same MSE in the first some iterations, then gradually the MSE of the Smart Gradient decreases showing a clear improvement, and it stays lower than the MSE of the Vanilla Gradient till the end of the iterations. Matrix $\pmb {G}^{(k)}$ needs at least $n$ iterations to be filled properly with informative directions, and improvements of the rounding off errors become evident after the number of iterations exceeds the dimension of $\pmb x$. This matrix continues rolling up for different $k$ until the optimum is found.  

It is clear that the gradient is calculated more accurately by Smart Gradient for all the functions in this optimization framework, as the prior information we have about the the descent directions boosts this accuracy. We examined the use of Smart Gradients on range of different test
 functions, all showing similar behaviour and overall improvement to
 the examples reported. For higher order finite different schemes with
 increased accuracy, the Smart Gradients had similar overall MSE than
 estimating Vanilla gradients directly. This is reasonable, as Smart
 Gradients does not offer any improvement in the limit.

\subsection{Smart Hessian}

The Smart Gradient technique can be extended and applied on Hessians. Assume for a continuous function $f(\pmb x)$, the second partial derivatives exit, we can estimate the Hessian $\nabla^2 f^{(k)}(\pmb x^{(k)})$ of this objection function at iteration $k$ based on some directions $\{{\Tilde{\pmb d}}_{1}^{(k)} | {\Tilde{\pmb d}}_{2}^{(k)} | \ldots | {\Tilde{\pmb d}}_{n}^{(k)} \}$ as we did for gradient in formulas \eqref{gradform4} and \eqref{gradform5},
\begin{equation}
 	\Tilde{\nabla}^2_{\Tilde{\pmb d}} f^{(k)}(\pmb x^{(k)}) = {{\pmb {\Tilde{G}}^{(k)}}}^{-T} \Tilde{\nabla}^2 h^{(k)}(\pmb \varphi)\Big|_{\pmb \varphi = 0}  {{\pmb {\Tilde{G}}^{(k)}}}^{T} \text{ where } h(\pmb \varphi) = f(\pmb  x^{(k)} + \pmb {\Tilde{G}}^{(k)} \pmb \varphi)
 	\label{hessformula}
 \end{equation}

One important application to this Smart Hessian is computing the Hessian of a function at its mode. This adaptive technique can help describing better the curvature of a function at its optimum, using the last $n$ descent directions.

\subsection{Autoregressive time-series model with \texttt{R-INLA}}\label{timeseries}

A motivation application to the Smart Gradient technique is fitting a model in Bayesian Inference using Integrated Nested Laplace Approximation (INLA) method \cite{Rue2009ApproximateBI}, \cite{Rue2016BayesianCW}. INLA uses combinations of analytical approximations and numerical integration to obtain approximated posterior distributions of the paramaters. It uses the BFGS method \cite{Fletcher1988PracticalMO} to reach the optimum value of the hyperparamater vector $\pmb \theta$ in the model. At each iteration, an inner optimization takes place to approximate the posterior distribution of the latent field $\pmb x$ by Gaussian approximation and get the mode $\pmb x^*$. This prohibits exact function evaluations.

The proposed technique of Smart Gradient is already used as the default option in \texttt{INLA} package in \texttt{R} to optimize the hyperparameters, using the following argument \texttt{inla(control.inla=list(use.directions = TRUE))} in \texttt{inla} function. We show in the next example how this method is used to fit a time series model using the INLA method.

\begin{figure}[hbt!]
\centering
    \includegraphics[scale=0.3]{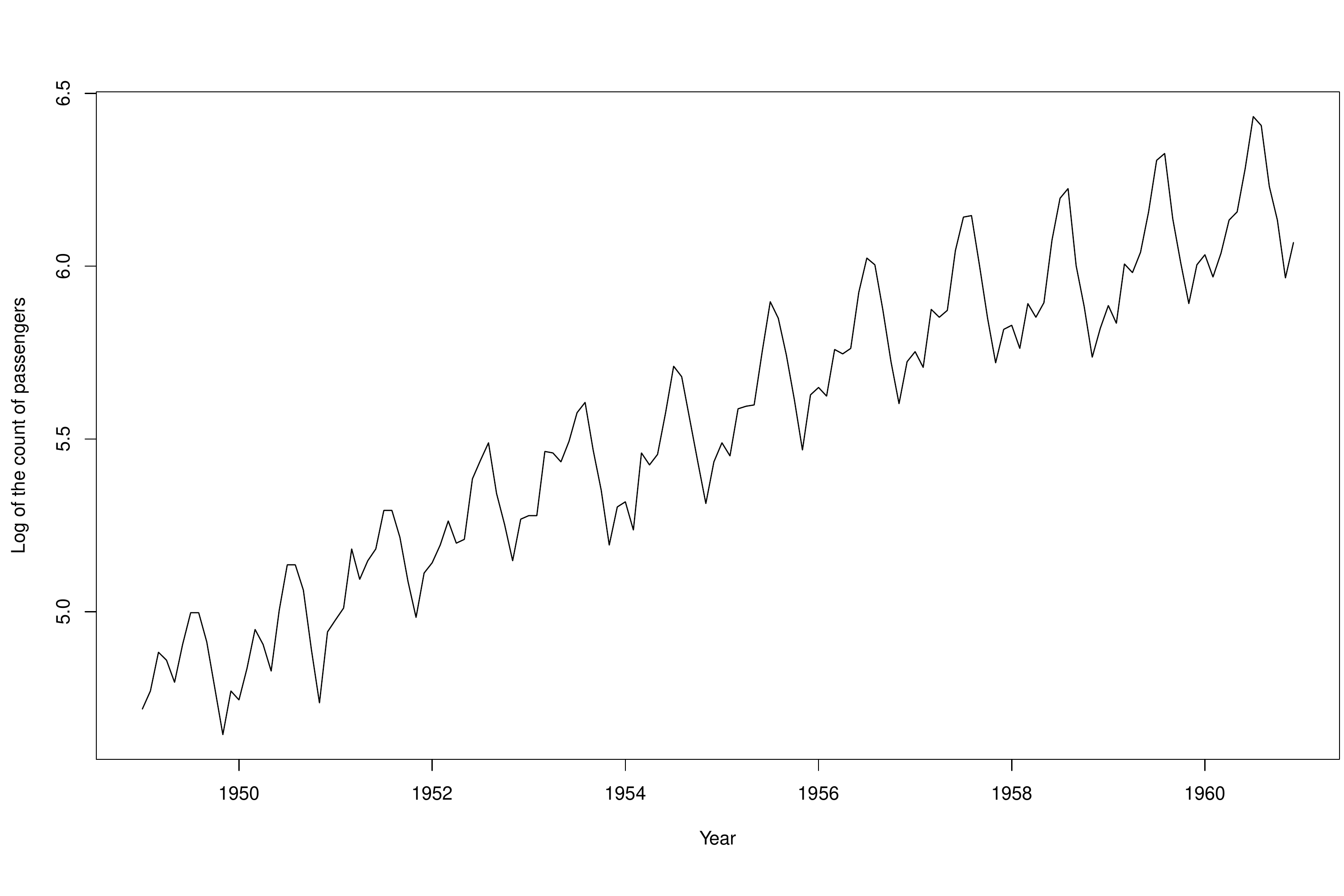}
    \caption{The log of air passengers from 1949 to 1961}
    \label{fig:airpass}
\end{figure}

Consider the R dataset of monthly totals of international airlines passengers between 1949 to 1960, see Figure \ref{fig:airpass}, one possible option is to fit the response as autoregressive model of order 1,

\begin{center}
$\log(\pmb y) \sim \beta \pmb t + \pmb u + \pmb \epsilon$
\end{center}

\begin{center}
    $u_i \sim \phi u_{i-1} + \varepsilon_i$, $\varepsilon_i \sim \mathcal{N}(0, \tau_u^{-1})$, $i = 2,\ldots, n$, and \\$u_1 \sim \mathcal{N}(0, (\tau_u(1 - \phi^2))^{-1})$
\end{center}
where $\pmb y$ is the number of air passengers, $\pmb t$ is the year, $\epsilon_j$ follows skew-normal distribution with zero mean, $j = 1,\ldots, n$, standardised skewness $\gamma$ and $\tau_y$ as the precision parameter. The fixed effect $\beta$ on $\pmb t$ follow a weakly informative Gaussian distribution a priori, with constant $\tau_\beta$ as precision. The four hyperparameters $(\tau_y, \gamma, \tau_u, \phi )$ in our model are transformed to $\pmb \theta = (\theta_1 = \log(\tau_y),\theta_2 = \log((1 + \gamma / 0.988) / (1 - \gamma / 0.988)), \theta_3 = \log(\tau_u(1 - \phi^2)), \theta_4 = \log((1 + \phi)/(1 - \phi)))$ for unconstrained optimization.

 In this inferential procedure we need to find the mode of the hyperparameters $\pmb \theta$, so that the marginals of the elements of the Gaussian latent field $\pmb x^T = (\beta, \pmb u^T)$ can be estimated.  The objective function $f(\pmb \theta|\pmb y)$,
 $$ \log f(\pmb \theta|\pmb y) \approx - \log\pi(\pmb y|\pmb x^*,\pmb \theta) + \log\pi(\pmb x^*|\pmb \theta) + \log\pi(\pmb \theta) - \pi(\pmb x^*|\pmb y,\pmb \theta)$$

\noindent where $\pmb x^*$ is the mode we get from the inner optimization of $\pi(\pmb x|\pmb y,\pmb \theta)$ as it is approximated by Gaussian. Using PC priors \cite{Simpson2014PenalisingMC}, we take the prior distributions for  $\pmb \theta$, 
\begin{equation}
    \begin{split}
        \log \pi(\pmb \theta) &\approx -\log(0.2) e^{-\theta_1/2} -\theta_1/2
        \\& ~~~ + \log(\phi(\theta_2)) + \log(\Phi(10^{1/3}\theta_2)) 
        \\& ~~~ + \theta_3  -(5\times 10^{-5})e^{\theta_3}-0.075\theta_4^2
    \end{split}
\end{equation}
\noindent where $\phi(.)$ is the standard normal distribution and $\Phi(.)$ is the standard cumulative Gaussian distribution. With $\pmb \eta$ as the linear predictor and $\pmb Q^{*} = \pmb Q(\pmb \theta) + \pmb Q_l$, where $\pmb Q(\pmb \theta)$ is a combination of the AR1 model precision matrix $\pmb R(\pmb \theta)$ and $\tau_\beta$, 

$$ \pmb Q(\pmb \theta) = \begin{pmatrix}
\tau_{\beta} & 0\\
0 & \pmb R(\pmb \theta) \\
\end{pmatrix},$$
\noindent and $\pmb Q_l$ is the precision we get from the model, then
$$\pi(\pmb y|\pmb x,\pmb \theta) \approx  \exp\Big(-\displaystyle\frac{\tau_\epsilon}{2}(\pmb y - \pmb \eta)^{T}(\pmb y - \pmb \eta)\Big)$$
$$\pi(\pmb x|\pmb \theta) \propto \exp\Big(-\displaystyle\frac{1}{2}\pmb x^{T} \pmb Q(\pmb \theta)\pmb x\Big) $$
$$\pi(\pmb x|\pmb \theta,\pmb y) \approx  \exp\Big(-\displaystyle\frac{1}{2}(\pmb x - \pmb x^*)^{T} \pmb Q^{*}(\pmb x - \pmb x^*)\Big)$$

Due to the intractable form of the objective function $f(\pmb \theta|y)$, it is hard to get the exact gradient, so it needs to be estimated numerically. Here we use Smart Gradient integrated with the central difference method which is more accurate compared to using Vanilla Gradient with canonical basis. 

Using Smart Gradient technique in this unconstrained optimization, we get $\pmb \theta^* = (2454.55, -0.001, 43.483 , 0.744)$, and the Smart Hessian of the objective function $f$ at $\pmb \theta^*$ can be calculated easily the same way.
 
 \subsection{Continuous spatial statistical model with the stochastic partial differential equation (SPDE) method \texttt{R-INLA}}
\begin{figure}[hbt!]
\centering
    \includegraphics[scale=0.5]{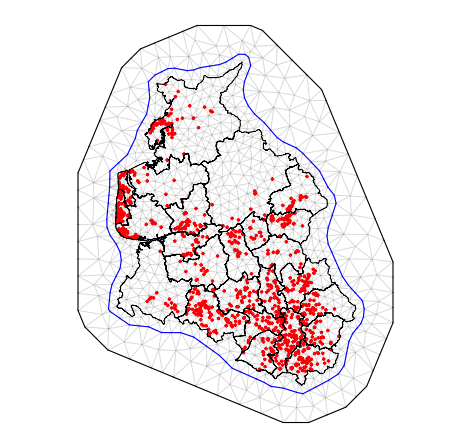}
    \caption{Exact residential locations of patients with AML}
    \label{fig:amlfig}
\end{figure}

Consider the R dataset \texttt{Leuk} that features the survival times of patients with acute myeloid leukemia (AML) in  Northwest England between 1982 to 1998. Exact residential locations and districts of the patients are known and indicated by the dots in Figure \ref{fig:amlfig}. 
The aim is to model the survival time based on various covariates $\pmb X$ and space $\pmb s$, with a Weibull model with shape $\alpha$,
\begin{center}
$t(\pmb s) \sim \text{Weibull} (\pmb\lambda = \exp[\pmb\beta\pmb X + \pmb u(\pmb s)], \alpha).$
\end{center}

The shape parameter of the likelihood is also considered a hyperparameter (different to Section \ref{timeseries}). To include a spatial element we use a Gaussian random effect $
\pmb u$ with a Matern covariance structure that has marginal variance $\sigma^2_u$ and nominal range $r = 2/\kappa$ such that the entry of the precision matrix of $\pmb u$ between locations $s_i$ and $s_j$ is $\pmb R_{ij} = \sigma^2_u e^{-\kappa (s_i-s_j)^{T}(s_i-s_j)}$. 

We use the finite element method and the mesh presented in Figure \ref{fig:amlfig} to estimate this model (for more details regarding the SPDE approach to continuous spatial modeling see \cite{krainski2018advanced}). In this example we have incomplete observations in the sense that some patients are still alive at the end of the study and some times are thus censored, indicated in the variable $\pmb c$. The data is thus $\pmb y = \{\pmb t, \pmb c\}$ instead of a univariate observation as in Section \ref{timeseries}.

Thus we have three hyperparameters in this model $\pmb\theta = \{10\log(\alpha),$\\ $\log(2/\kappa), \log(\sigma^{-2}_u)\}$, one from the likelihood and two from the spatial field which we need to optimize in order to obtain the marginal posteriors of the latent field. The objective function for $\pmb\theta$ is obtained similar to that in Section \ref{timeseries} and again is intractable, necessitating a numerical gradient. We deploy Smart Gradient technique again and we calculate the optimal values for $\pmb\theta$ as $\pmb \theta^* = (-5.2246, 0.6530 , 0.1599)$.

\section{Discussion and further considerations}\label{Discussion_section}

Gradients are important tool in various applied mathematical and statistical methods. This wide use of gradients necessitates the search for techniques that improve its estimation. Here, we presented a simple framework that enhances the accuracy of gradient estimation within optimization context using the most recent differences between $\pmb x$ positions as directions to compute the gradient.

In general, using more precise gradient will lead to a better performance for gradient based optimization methods. For a moderate number of dimension, Smart Gradient method shows improvement in its accuracy with essentially no cost. For instance, a Smart Gradient technique with first order central difference method, where function is evaluated at only two evaluation points can perform as well as higher order differences methods which are computationally more costly. 

The proposed method is naturally extended to improve the estimate of the Hessian through the attained descent directions. It is implemented in the accompanying R package \textbf{\texttt{smartGrad}} available on the github at \url{esmail-abdulfattah/Smart-Gradient}. Additionally, \textbf{\texttt{smartGrad}} can be used to enhance a user-defined gradient formula through the function \textbf{\texttt{makeSmart}}, see Appendix \ref{AppendixA}. C\texttt{++} code is also available on github. 


\vspace{0.5cm}


\newpage
\bibliography{mybib}{}
\bibliographystyle{elsarticle-num} 

\newpage
\appendix
\section{smartGrad Installation and Examples} \label{AppendixA}

\subsection{Installation}
\begin{verbatim}
    library("devtools")
    install_github("esmail-abdulfattah/Smart-Gradient", 
                    subdir = "smartGrad")
\end{verbatim}

\subsection{makeSmart Function}

To illustrate how to use this function, we use the Extended Rosenbrock function of dimension $n$. It has global minimum at $\pmb x^* = \pmb 1$. We estimate the gradient of this objective function $f$ using a simple central difference method, with step size $10^{-3}$.

\vspace{0.4cm}
\begin{verbatim}
    myfun <- function(x) {
      res <- 0.0
      for(i in 1:(length(x)-1))
        res <- res + 100*(x[i+1] - x[i]^2)^2 + (1-x[i])^2
      return(res)
    }

    mygrad <- function(fun,x){
      h = 1e-3
      grad <- numeric(length(x))
      for(i in 1:length(x)){
        e = numeric(length(x))
        e[i] = 1
        grad[i] <- (fun(x+h*e) - fun(x-h*e))/(2*h)
      }
      return(grad)
    }
    
\end{verbatim}

\noindent A user can change a numerical gradient function \texttt{mygrad} to a SMART numerical gradient function \texttt{mySmartgrad}, and then this SMART function can be used instead. We use the \texttt{optim} function from \texttt{stats} package with \texttt{BFGS} algorithm. 

\vspace{0.4cm}
\begin{verbatim}
    library("stats")
    library("smartGrad")
    x_dimension = 5
    x_initial = rnorm(x_dimension)
    result <- optim(par = x_initial, 
                            fn = myfun,
                            gr = makeSmart(fn = myfun,gr = mygrad),
                            method = c("BFGS"))
\end{verbatim}

\newpage

\end{document}